\documentclass{proc} 

\usepackage[a5paper,margin=0.5in,hmargin={0.8in,0.5in},includehead]{geometry}

\usepackage{multicol}
\usepackage{amsmath}

\usepackage{amssymb}
\usepackage{wrapfig}
\usepackage{caption,setspace,subfigure}

\let\footnote\savefootnote
\let\footnotetext\savefootnotetext
\let\footnoterule\savefootnoterule
\normallatexbib 

\usepackage{graphicx}
\usepackage{edbkps} 

\begin{document}

\articletitle[Strengthening SONC Relaxations]
{Strengthening SONC Relaxations\\
with Constraints Derived from \\
Variable Bounds\thanks{The work for this article has been conducted within the Research Campus
Modal funded by the German Federal Ministry of Education and Research (BMBF grant numbers 05M14ZAM,
05M20ZBM).}
}

\author{Ksenia Bestuzheva,\altaffilmark{1} Ambros Gleixner,\altaffilmark{1,2}
and Helena V{\"o}lker \altaffilmark{1}}

\altaffiltext{1}{Zuse Institute Berlin,
Department AIS$^2$T,
Takustr. 7, 14195 Berlin,
Germany}

\altaffiltext{2}{HTW Berlin,
Treskowallee 8, 10318 Berlin,
Germany}


\anxx{Ksenia Bestuzheva, Zuse Institute Berlin\,
Takustra{\ss}e 7\, 14195 Berlin\,
Germany\, \emailfont{bestuzheva@zib.de}}

\anxx{Ambros Gleixner, HTW Berlin\,
Treskowallee 8\, 10318 Berlin\,
Germany\,
\emailfont{gleixner@zib.de}}

\anxx{Helena V{\"o}lker, Zuse Institute Berlin\,
Takustra{\ss}e 7\, 14195 Berlin\,
Germany\, \emailfont{voelker@zib.de}}

\begin{abstract}
Certificates of polynomial nonnegativity can be used to obtain tight dual bounds for polynomial
optimization problems.
We consider Sums of Nonnegative Circuit (SONC) polynomials certificates, which are well
suited for sparse problems since the computational cost depends only on the number of terms in the
polynomials and does not depend on the degrees of the polynomials.
This work is a first step to integrating SONC-based relaxations of polynomial problems into a
branch-and-bound algorithm.
To this end, the SONC relaxation for constrained optimization problems is extended in order to
better utilize variable bounds, since this property is key for the success of a relaxation in the
context of branch-and-bound.
Computational experiments show that the proposed extension is crucial for making the SONC
relaxations applicable to most constrained polynomial optimization problems and for integrating the
two approaches.
\end{abstract}


\section{Introduction}

Consider a constrained polynomial optimization problem of the form
\begin{subequations}
\label{eqn:optproblem}
\begin{align}
\min_{\mathbf{x} \in \mathbb{R}^n}~~   &f(\mathbf{x}) = \sum_{\alpha \in \mathcal{A}_0} f_{\alpha} \mathbf{x}^{\alpha} \label{eqn:optproblemObj} \\
\text{s.t.}~~ &g_i(\mathbf{x}) = \sum_{\alpha \in \mathcal{A}_i} g_{i,\alpha} \mathbf{x}^{\alpha} \geq 0, \hspace{5mm} i=1,\ldots,m, \label{eqn:optproblemConst}
\end{align}
\end{subequations}
where $f_{\alpha}, g_{i,\alpha} \in \mathbb{R}$ are nonzero coefficients of the polynomials, and
$\mathcal{A}_i \subset \mathbb{N}^n$, $i = 0,\dots,m$, are supports of polynomials $f$ and $g_i$.

A monomial is written as $\mathbf{x}^{\alpha}:=
x_{1}^{\alpha_{1}}\cdot \ldots \cdot x_{n}^{\alpha_{n}}$.
Exponent $\alpha$ is referred to as a \emph{monomial square} of $f$ if the term $f_{\alpha}\mathbf{x}^{\alpha}$ satisfies $f_{\alpha}\geq 0$ and $\alpha \in (2\mathbb{N})^n$.
The \emph{Newton polytope} New$(f)$ of a polynomial $f$ with support $\mathcal{A}$ is defined to be the convex hull of the exponents of $f$, i.e.
$\text{New}(f) = \text{conv}(\mathcal{A})$.
Let $\Delta(\mathcal{A}) = \mathcal{A} \setminus \text{New}(f)$.

\begin{definition} \label{def:STpoly}
Consider a polynomial $f$ with support $\mathcal{A}$ and let $k = |\text{New}(f)|$.
$f$ is an \emph{ST-polynomial}~\cite{Lower} if it has the form
\begin{equation} \label{eqn:STpoly}
f(\textbf{x}) = \sum_{j=0}^k f_{\alpha(j)} x^{\alpha(j)} + \sum_{\beta \in \Delta(\mathcal{A})} f_{\beta} x^{\beta}
\end{equation}
such that New$(f)$ is a simplex whose vertices $\alpha(j)$, $j = 0,\dots,k$, are monomial squares,
 and for all $\beta \in \Delta(\mathcal{A})$ there exist $\lambda_0^{(\beta)},\ldots,\lambda_k^{(\beta)} \geq 0$ forming the unique convex combination
\begin{equation} \label{eqn:barycoord}
\sum_{j=0}^k \lambda_j^{(\beta)} = 1 \hspace{5mm} \text{ and } \hspace{5mm} \sum_{j=0}^{k} \lambda_j^{(\beta)} \alpha(j) = \beta.
\end{equation} 

The values $\lambda_0^{(\beta)},\ldots,\lambda_{k}^{(\beta)} \in \mathbb{R}$ are referred to as \textit{barycentric coordinates} of $\beta$,
$\beta \in \Delta(A)$ are referred to as \emph{inner terms} of $f$,
and the vertex set $\{\alpha(0),\ldots,\alpha(k)\}$ of the simplex is referred to as a \emph{cover} of the inner term $\beta$.
We will say that $\beta$ is \emph{covered} by $\{\alpha(0),\ldots,\alpha(k)\}$.
\end{definition}

A special case of ST-polynomials are circuit polynomials~\cite{SONC}, which are defined
as ST-polynomials with one inner term $\beta$.
For circuit polynomials, the barycentric coordinates are used to define the circuit number of an inner term $\beta$.
Nonnegativity can then be decided by comparing the circuit number to the coefficient of the inner term.

Optimization problems can be stated as nonnegativity problems, since one can write $\min f(x)$ s.t.
$g(x)\geq 0$ equivalently as $\sup \{ \gamma \in \mathbb{R}: \hspace{2mm} f(x)-\gamma\geq 0, g(x)
\geq 0 \}$.
One way of certifying nonnegativity of a polynomial is to write it as a sum of nonnegative circuit
polynomials.

Using geometric programming, one can obtain a lower found on the optimal value of problem~\eqref{eqn:optproblem} by finding a
SONC decomposition of the Lagrangian of problem~\eqref{eqn:optproblem}.

\section{Enhanced SONC Relaxations}

A SONC decomposition can only be found if all vertices of New$(f)$ are monomial squares.
Therefore, when working with an arbitrary polynomial, it may be necessary to reformulate the lower
bounding problem so that each term that is not a monomial square is covered by even exponents.
Moreover, a desired property in a relaxation that is solved as part of a branch-and-bound algorithm
is its ability to utilise tighter variable bounds in order to obtain lower bounds of improved quality.

We assume that finite lower and upper bounds are provided for each variable, i.e.
$\max \{\vert l_i \vert, \vert u_i \vert \} < \infty$ for each $i=1,\ldots,n$.
The bounds are utilized in order to derive polynomial constraints of the form
\begin{equation*} \label{eqn:boundcons}
x_{i}^{a_i} \leq \max \{\vert l_i \vert, \vert u_i \vert \}^{a_i},
\end{equation*} 
referred to as \emph{bound constraints}.

Let the monomial squares of the Lagrangian of problem~\eqref{eqn:optproblem} be denoted by $\alpha(0),\ldots,\alpha(k)$ and assume w.l.o.g. $\alpha(0) = 0$.
We require a set of exponents $a_i$ such that these exponents together with the even exponents $\alpha(0),\ldots,\alpha(k)$ give a valid cover for each $\beta\in\Delta(A)$:

\begin{align*}\label{eqn:OPexponent}
&\left\{ a_1,\ldots,a_n \mid \forall~\beta \in \Delta(A)~\exists~\lambda_0^{(\beta)},\ldots,\lambda_{k+n}^{(\beta)}\geq 0: \vphantom{\sum_{i=1}^{n+k}} \right. \\ 
&\hspace{4mm}\left. a_i \in 2 \mathbb{N},~\beta= \sum_{j=0}^{k} \lambda_j^{(\beta)} \alpha(j) + \sum_{i=1}^{n} \lambda_{k+i}^{(\beta)} \mathbf{a^i} \text{ and } \sum_{j=0}^{k+n} \lambda_j^{(\beta)}= 1 \right\}. \nonumber
\end{align*}
We define $\mathbf{a^i}$ as the vector with $0$ in all components except the i-th, where it is $a_i$.
It is easy to check that this vector satisfies the above condition.

\section{Computational Results}

Table~\ref{tab:boundsYesNo} compares the results produced when using the relaxation enhanced with
bound constraints with the results produced when using the standard SONC relaxation.
\begin{center}
\begin{table}[ht]
\begin{tabular*}{\linewidth}{@{\extracolsep{\fill}}ccccc}
\sphline
\it Test Run &  \multicolumn{4}{c}{\it Solution Status of the Relaxation}  \\ 
& \it optimal & \it infeasible & \it unsolvable & \it numerical error   \\
\sphline
With BndConss & 330 & 0 & 0 & 19 \\
Without BndConss & 9 & 4 & 329 & 7 \\
\sphline
\end{tabular*}
\caption[Comparison of the Relaxator with and without Bounds]{Comparison of the relaxator with bound constraints and without bound constraints.} \label{tab:boundsYesNo}
\end{table}
\end{center}
We can see that for the majority of instances, the standard SONC relaxation fails to find a lower bound.
This is due to the fact that for most instances the Lagrangian does not have a Newton polytope with even vertices. 
Adding the bound constraints increases the number of instances where a finite lower bound was found from $2.6\%$ to $94.6\%$, as the bound constraints always provide a valid cover.
Thus, our approach drastically improves the applicability of SONC relaxations to constrained polynomial optimization problems.

\bibliographystyle{plain}
\chapbblname{SONC_HUGO2022}  
\chapbibliography{SONC_HUGO2022}

\end{document}